\documentclass[12pt]{article}
\usepackage{a4, amssymb, amsmath, amsthm, xypic}

\begin{document}

\title{Classical dynamical r-matrices,
Poisson homogeneous spaces, and Lagrangian subalgebras}

\author{Eugene~Karolinsky\footnote{Research was
supported in part by CRDF grant UM1--2091 and the Royal
Swedish Academy of Sciences.}\\ {\small \tt
Department of Mathematics, Kharkov National University,}\\ {\small
\tt 4 Svobody Sq., Kharkov, 61077, Ukraine}\\ {\small \tt
karol@skynet.kharkov.com}\\  \\
Alexander~Stolin\\ {\small \tt Department of Mathematics,
University of G\"oteborg,}\\ {\small \tt SE-412 96 G\"oteborg,
Sweden}\\ {\small \tt astolin@math.chalmers.se}}

\date{}


\newcommand{\g}{\mathfrak g}
\newcommand{\h}{\mathfrak h}
\renewcommand{\l}{\mathfrak l}
\newcommand{\m}{\mathfrak m}
\newcommand{\n}{\mathfrak n}
\renewcommand{\u}{\mathfrak u}
\newcommand{\z}{\mathfrak z}
\newcommand{\p}{\mathfrak p}
\renewcommand{\v}{\mathfrak v}
\newcommand{\Dg}{D(\g)}
\newcommand{\geps}{\g [\varepsilon]}
\newcommand{\Ker}{\mbox{\rm Ker} \, }
\newcommand{\Lie}{\mbox{\rm Lie} \, }
\newcommand{\CYB}{\mbox{\rm CYB}}
\newcommand{\Alt}{\mbox{\rm Alt}}
\newcommand{\Ad}{\mbox{\rm Ad}}
\newcommand{\ad}{\mbox{\rm ad}}
\newcommand{\id}{\mbox{\rm id}}
\newcommand{\Span}{\mbox{\rm Span}}
\newcommand{\gxg}{\g \! \times \! \g}
\newcommand{\gotg}{\g \otimes \g}
\newcommand{\gwg}{\wedge^2 \g}
\newcommand{\gotgu}{(\g \otimes \g)^\u}
\newcommand{\mwm}{\wedge^2 \m}
\newcommand{\mwmu}{(\wedge^2 \m)^\u}
\newcommand{\scpr}{\langle \cdot \, , \cdot \rangle}
\newcommand{\RRR}{{\bf R}}
\newcommand{\NNN}{{\bf N}}
\newcommand{\UUU}{{\bf U}}
\newcommand{\Pairs}{\mbox{\bf Pairs}(\g)}
\newcommand{\Lagr}{\mbox{\bf Lagr}(\g)}
\newcommand{\Pairsu}{\mbox{\bf Pairs}(\g, \u)}
\newcommand{\Lagru}{\mbox{\bf Lagr}(\g, \u)}
\newcommand{\Dynr}{\mbox{\bf Dynr}(\g, \u, \Omega)}
\newcommand{\Modr}{{\cal M}(\g, \u, \Omega)}
\newcommand{\Momega}{{\cal M}_\Omega}
\newcommand{\Maps}{\mbox{\bf Map}(D, G)^\u}
\newcommand{\Mapszero}{\mbox{\bf Map}_0(D, G)^\u}
\newcommand{\Dom}{\mbox{\rm Domain}}
\newcommand{\halfomega}{\frac{\Omega}{2}}
\newcommand{\modomega}{\frac{1}{4}[\Omega^{12}, \Omega^{23}]}
\newcommand{\nh}{\NNN,\,h}

\newtheorem{proposition}{Proposition}
\newtheorem{theorem}[proposition]{Theorem}
\newtheorem{corollary}[proposition]{Corollary}
\newtheorem{lemma}[proposition]{Lemma}
\newtheorem{remark}[proposition]{Remark}
\newtheorem{example}[proposition]{Example}


\maketitle

\begin{abstract}
In \cite{Lu CDYBE} Lu showed that any dynamical r-matrix for the
pair $(\g,\u)$ naturally induces a Poisson homogeneous structure
on $G/U$. She also proved that if $\g$ is complex simple, $\u$ is
its Cartan subalgebra and $r$ is quasitriangular, then this
correspondence is in fact 1-1. In the present paper we find some
general conditions under which the Lu correspondence is 1-1. Then
we apply this result to describe all triangular Poisson
homogeneous structures on $G/U$ for a simple complex group $G$ and
its reductive subgroup $U$ containing a Cartan subgroup.
\end{abstract}

\begin{sloppy}

\section{Introduction}

The notion of a Poisson-Lie group was introduced almost 20 years
ago by Drinfeld in \cite{Drinfeld PL}.  Its infinitesimal
counterpart, Lie bialgebras, were introduced in the same paper and
later it was explained that these objects are in fact
quasiclassical limits of quantum groups (see \cite{Drinfeld QG}).
Lie bialgebra structures on a Lie algebra $\g$ are in a natural
1-1 correspondence with Lie algebra structures on the vector space
$\Dg = \g \oplus \g^*$ with some compatibility conditions. $\Dg$
with this Lie algebra structure is called the double of the Lie
bialgebra $\g$.

The most popular and important class of Lie bialgebras is the
class of quasitriangular Lie bialgebras (see \cite{Drinfeld QG}).
They can be defined by an element $r \in \g \otimes \g$ (called
the classical r-matrix) such that
$$\Omega := r + r^{21}$$
is $\g$-invariant, and the classical Yang-Baxter equation (CYBE)
$$[r^{12},r^{13}]+[r^{12},r^{23}]+[r^{13},r^{23}] = 0$$
is satisfied. If $r$ is skew-symmetric, then one says that the
corresponding Lie bialgebra is triangular. In general, $\Lambda :=
r - \halfomega$ (i.e., the skew-symmetric part of $r$) satisfies
the modified CYBE
$$[\Lambda^{12}, \Lambda^{13}] + [\Lambda^{12}, \Lambda^{23}] +
[\Lambda^{13}, \Lambda^{23}] = \modomega.$$

It is well known (and can be easily shown) that if $\g$ is a
complex simple finite-dimensional Lie algebra, then any Lie
bialgebra structure on $\g$ is quasitriangular. For the case
$\Omega \neq 0$ (``quasitriangular case in the strict sense") they
were classified by Belavin and Drinfeld, see \cite{Bel-Dr FA,
Bel-Dr}. The triangular case was studied in \cite{Stolin0,
Stolin2, Stolin3}.

In the paper \cite{Stolin1} it was shown that for such $\g$ there
are only two possible structures of the $\Dg$. In the triangular
case $\Dg = \geps = \g \oplus \g\varepsilon$, where $\varepsilon^2
= 0$ and otherwise, $\Dg = \gxg$ (and $\g$ is embedded diagonally
into $\gxg$). Then it is clear that skew-symmetric solutions of
the CYBE (resp.\ the modified CYBE with $\Omega \neq 0$) are in a 1-1
correspondence with Lagrangian subalgebras $\l$ in $\geps$ (resp.\
in $\gxg$) such that $\l \cap \g = 0$.

Along with the Poisson-Lie groups it is natural to study their
Poisson actions. The notion of Poisson action was introduced by
Semenov-Tian-Shansky in \cite{STS}. Poisson homogeneous
spaces are of our special interest. Drinfeld in \cite{Drinfeld Lagr}
gave a general approach
to the classification of Poisson homogeneous spaces. Namely, he
showed that if $G$ is a Poisson-Lie group, $\g$ is the
corresponding Lie bialgebra, then Poisson homogeneous $G$-spaces
are essentially in a 1-1 correspondence with $G$-orbits on the set
of all Lagrangian subalgebras in $\Dg$. A classification of
Lagrangian subalgebras in some important cases (including the case
$\g$ is complex simple, $\Dg = \gxg$) was obtained in \cite{Kar2,
Kar3, Kar1}.

At the same time an important generalization of the CYBE, the dynamical
classical Yang-Baxter equation, was introduced in physics and mathematics.
Notice that this equation is defined for a pair $(\g,\u)$, where $\u$ is a
subalgebra of $\g$. From the mathematical point of view it was presented by
Felder in \cite{Felder, Felder1}. This equation and its quantum analogue
were studied in many papers, see \cite{E-V lectures, E-V, Sch, Xu, Xu1}. First
classification results for the solutions of the classical dynamical
Yang-Baxter equation (dynamical r-matrices) were obtained by Etingof,
Varchenko, and Schiffmann in \cite{E-V, Sch}.

Later Lu \cite{Lu CDYBE} found a connection (which is essentially a 1-1
correspondence) between dynamical r-matrices for the pair $(\g, \u)$ (where
$\u $ is a Cartan subalgebra of the complex simple finite-dimensional
algebra $\g$), and Poisson homogeneous $G$-structures on $G/U$. Here $U
\subset G$ are connected Lie groups corresponding to $\u \subset \g$, and
$G$ is equipped with the standard quasitriangular (with $\Omega \neq 0$)
Poisson-Lie structure.

Lu also noticed that this connection can be generalized to the case $\u$ is
a subspace in a Cartan subalgebra (with some ``regularity'' condition). The
dynamical r-matrices for the latter case were classified by Schiffmann in
\cite{Sch}. In this case connections between dynamical r-matrices and
certain Lagrangian subalgebras can be derived directly from \cite{Sch}.

Now let $G$ be a connected complex  semisimple Lie group, and let $U$ be its
connected subgroup. Suppose $\u \subset \g$ are the corresponding Lie
algebras. In the present paper we consider connections between Poisson
homogeneous structures on $G/U$ related to the triangular Poisson-Lie
structures on $G$ (i.e., with $\Omega=0$), where $U$ is a reductive subgroup
containing a Cartan subgroup of $G$, and triangular dynamical r-matrices for
the pair $(\g,\u)$.

In fact, our results are based on a general result on relations between
classical dynamical r-matrices and Poisson homogeneous structures (see
Theorem \ref{twistdyn1}), which is valid also in the quasitriangular case.
Notice that the results of Sections \ref{sectdyn} and \ref{secttwist} can be
used to describe a 1-1 correspondence between dynamical r-matrices for the
pair $(\g, \u)$, where $\u \subset \g$ is a Cartan subalgebra, and Poisson
homogeneous $G$-structures on $G/U$, where $G$ is equipped with any
quasitriangular (with $\Omega \neq 0$) Poisson-Lie structure (of course the
latter result is due to Lu). Our approach is based on some strong
classification results for dynamical r-matrices given recently by Etingof
and Schiffmann in \cite{E-S}.

The paper is organized as follows. In Section \ref{sectdyn} we describe a
correspondence between the (moduli space of) dynamical r-matrices for a pair
$(\g, \u)$ and Poisson homogeneous $G$-structures on $G/U$ proving that
under certain assumptions it is a bijection. In Section \ref{secttwist} we
consider a procedure of twisting for Lie bialgebras and examine its impact
on the double $D(\g)$ and Poisson homogeneous spaces for corresponding
Poisson-Lie groups. Then we use the twisting to weaken some restrictions
needed in Section \ref{sectdyn}. In Section \ref{sectex} we consider the
basic example of our paper: $\g$ is semisimple, $\u \subset \g$ is a
reductive Lie subalgebra that contains some Cartan subalgebra of $\g$, and
the Lie bialgebra structure on $\g$ is triangular (i.e., $\Dg = \geps$).
In Appendix A we present a general approach to the description of all
Lagrangian subalgebras in $\geps$ and give a direct classification of the
Lagrangian subalgebras $\l \subset \geps$ such that $\l \cap \g = \u$.
In Appendix B we describe a method that allows to get constant r-matrices
from dynamical ones.

\subsection*{Acknowledgements} The authors are grateful to
Pavel Etingof and Olivier Schiffmann for valuable discussions.

\section{Classical dynamical r-matrices and Poisson homogeneous spaces}
\label{sectdyn}

In this section we assume $\g$ to be any finite-dimensional Lie algebra over
$\mathbb C$. Let $G$ be a connected Lie group such that $\Lie G = \g$. Let
$\u \subset \g$ be a Lie subalgebra (not necessary abelian). By $U$ denote
the connected subgroup in $G$ such that $\Lie U = \u$. We propose (under
certain conditions) a connection between dynamical r-matrices for the pair
$(\g, \u)$ and Poisson structures on $G/U$ that make $G/U$ a {\it Poisson}
homogeneous $G$-space (for certain Poisson-Lie structures on $G$). Note that
this connection was first introduced by Jiang-Hua Lu in \cite{Lu CDYBE} for
the case $\g$ is simple, $\u$ is a Cartan subalgebra, and the dynamical
r-matrix has non-zero coupling constant. Our result is inspired by \cite{Lu
CDYBE}.

In order to recall the definition of the classical dynamical
r-matrix we need some notation. Let $x_1,...,x_r$ be a basis of
$\u$. By $D$ denote the formal neighborhood of $0$ in $\u^*$. By
functions from $D$ to a vector space $V$ we mean the elements of
the space $V[[x_1,...,x_r]]$, where $x_i$ are regarded as
coordinates on $D$. If $\omega \in \Omega^k(D,V)$ is a $k$-form
with values in a vector space $V$, we denote by
$\overline{\omega}: D \to \wedge^k \u \otimes V$ the corresponding
function. Finally, for an element $r \in \gotg$ we define the
classical Yang-Baxter operator $$\CYB(r)=
[r^{12},r^{13}]+[r^{12},r^{23}]+[r^{13},r^{23}].$$

Recall that a {\it classical dynamical r-matrix} for the pair
$(\g, \u)$ is an $\u$-equivariant function $r: D \to \gotg$ that
satisfies the classical dynamical Yang-Baxter equation (CDYBE):
\begin{equation}
\Alt(\overline{dr}) + \CYB(r) = 0,
\end{equation}
where for $x \in \g^{\otimes 3}$, we set $\Alt(x) = x^{123} +
x^{231} + x^{312}$ (see \cite{E-V lectures, E-S, E-V}). Usually
one requires also an additional {\it quasi-unitarity condition}:
$$r + r^{21} = \Omega \in (S^2 \g)^\g.$$ Note that if $r$
satisfies the CDYBE and the quasi-unitarity condition, then
$\Omega$ is a constant function.

Suppose $\Omega \in (S^2 \g)^\g$. Let us denote by $\Dynr$ the set
of all classical dynamical r-matrices $r$ for the pair $(\g, \u)$
such that $r + r^{21} = \Omega$.

Denote by $\Maps$ the set of all regular $\u$-equivariant maps from $D$ to
$G$. Suppose $g \in \Maps$. For any $\mathfrak u$-equivariant function $r: D
\to \gotg$ set $$r^g = (\Ad_g \otimes \Ad_g)(r - \overline{\eta_g} +
\overline{\eta_g}^{21} + \tau_g),$$ where $\eta_g=g^{-1}dg$, and
$\tau_g(\lambda) = (\lambda \otimes 1 \otimes 1)([\overline{\eta_g}^{12},
\overline{\eta_g}^{13}](\lambda))$. Then $r$ is a classical dynamical
r-matrix iff $r^g$ is (see \cite{E-S}). The transformation $r \mapsto r^g$
is called a gauge transformation. It is indeed an action of the group
$\Maps$ on $\Dynr$ (i.e., $(r^{g_1})^{g_2} = r^{g_2g_1}$). Following
\cite{E-S} denote by $\Modr$ the moduli space $\Mapszero\backslash\Dynr$,
where $\Mapszero$ is the subgroup in $\Maps$ consisting of maps $g$
satisfying $g(0) = e$.

In what follows we need some notation. Suppose $a \in \mathfrak
g^{\otimes k}$. By $\overrightarrow{a}$ (resp.
$\overleftarrow{a}$) denote the left (resp. right) invariant
tensor field on $G$ corresponding to $a$.

Suppose $\rho \in \gotg$ satisfies the classical Yang-Baxter
equation (CYBE), i.e., $\CYB(\rho) = 0$. Assume also that $\rho +
\rho^{21} = \Omega$ (i.e., $\rho = \halfomega + \Lambda$, where
$\Lambda \in \gwg$). Introduce a bivector field $\pi_\rho =
\overrightarrow{\rho} - \overleftarrow{\rho} =
\overrightarrow{\Lambda} - \overleftarrow{\Lambda}$ on $G$. It is
well known that $(G, \pi_\rho)$ is a Poisson-Lie group.

Now let $r \in \Dynr$. We have $r = \halfomega + A$, where $A \in
\gwg$. Set $\tilde\pi_r := \overrightarrow{r(0)} -
\overleftarrow{\rho} = \overrightarrow{A(0)} -
\overleftarrow{\Lambda}$. Consider a bivector field $\pi_r$ on
$G/U$ defined by $\pi_r\left(\underline{g}\right) =
p_*\tilde\pi_r(g)$, where $p: G \to G/U$ is the natural
projection, and $\underline{g} = p(g)$. Note that $\pi_r$ is well
defined since $r(0) \in \gotgu$.

The following proposition belongs to Jiang-Hua Lu \cite{Lu CDYBE}
(note that in \cite{Lu CDYBE} it is stated for the case $\mathfrak
g$ is simple, $\u$ is a Cartan subalgebra, but the proof fits the
general case).

\begin{proposition} \label{Lu}
The bivector field $\pi_r$ is Poisson, and $(G/U, \pi_r)$ is a
Poisson homogeneous $(G, \pi_\rho)$-space. \hfill $\square$
\end{proposition}

\begin{proposition}
Suppose $g \in \Mapszero$. Then $\pi_r = \pi_{r^g}$.
\end{proposition}

\proof Since $(G/U, \pi_r)$ is a Poisson homogeneous $(G,
\pi_\rho)$-space, we see that $\pi_r$ depends only on $\pi_r
(\underline{e}) =$ the image of $r(0) - \rho$ in $\wedge^2 (\g /
\u)$. Thus it is enough to note that $r^g(0) - r(0) \in \u \otimes
\g + \g \otimes \u$. \hfill $\square$

\begin{corollary}
The correspondence $r \mapsto \pi_r$ defines a map from $\Modr$ to
the set of all Poisson $(G, \pi_\rho)$-homogeneous structures on
$G/U$. \hfill $\square$
\end{corollary}

Suppose now that the following conditions hold:

\begin{itemize}
\item[(a)] $\u$ has an $\u$-invariant complement
$\m$ in $\g$ (we fix one).
\item[(b)] $\Omega \in (\u \otimes \u) \oplus (\m
\otimes \m)$.
\item[(c)] $\rho \in \halfomega + \mwmu$.
\end{itemize}

\begin{theorem} \label{dyn}
Under the assumptions above the correspondence $r \mapsto \pi_r$
is a bijection between $\Modr$ and the set of all Poisson $(G,
\pi_\rho)$-homogeneous structures on $G/U$.
\end{theorem}

The rest of this section is devoted to the proof of Theorem
\ref{dyn}. First we recall some results from \cite{E-S}. Assume
that (a) holds. Set $$\Momega = \left\{ x \in \halfomega + \mwmu
\, \Bigm| \, \CYB(x) = 0 \ \mbox{in}\, \wedge^3(\g / \u)
\right\}.$$

\begin{theorem}[Etingof, Schiffmann \cite{E-S}] \label{ESmain}

1. Any class ${\cal C} \in \Modr$ has a representative $r \in \cal
C$ such that $r(0) \in \Momega$. Moreover, this defines an
embedding $\Modr \to \Momega$.

2. Assume that (b) holds. Then the map $\Modr \to \Momega$ defined
above is a bijection. \hfill $\square$
\end{theorem}

Now suppose $b \in (\wedge^2 (\g / \u))^\u = (\wedge^2 \mathfrak
m)^\u$. Set $\pi\left(\underline{g}\right) = (L_g)_*b +
p_*\pi_\rho(g)$. Since $\rho$ is $\u$-invariant, we see that
$\pi_\rho(g) = 0$ for $g \in U$; therefore $\pi$ is a well-defined
bivector field on $G/U$.

\begin{proposition} \label{poiss}
The bivector field $\pi$ is Poisson iff $\CYB(\rho + b) = 0$ in
$\wedge^3 (\g / \u)$.
\end{proposition}

\proof Set $a = \Lambda + b$. Define a bivector field
$\tilde{\pi}$ on $G$ by the formula $\tilde\pi =
\overrightarrow{a} - \overleftarrow{\Lambda}$. Note that
$\tilde\pi = \overrightarrow{b} + \pi_\rho$, therefore $\pi =
p_*\tilde\pi$. Let us normalize the Schouten bracket of the
bivector fields on $G$ in a way that $[\overrightarrow{x},
\overrightarrow{x}] = \overrightarrow{\CYB(x)}$ for all $x \in
\wedge^2 \g$. Then we have $$[\tilde\pi, \tilde\pi] =
[\overrightarrow{a}, \overrightarrow{a}] - 2[\overrightarrow{a},
\overleftarrow{\Lambda}] + [\overleftarrow{\Lambda},
\overleftarrow{\Lambda}] = \overrightarrow{\CYB(a)} -
\overleftarrow{\CYB(\Lambda)}.$$ Since $\rho = \halfomega +
\Lambda$ satisfies the CYBE, we see that $\CYB(\Lambda) = \modomega
\in (\wedge^3 \g)^\g$. Thus $$[\tilde\pi, \tilde\pi] =
\overrightarrow{\CYB(a) - \modomega} =
\overrightarrow{\CYB\left(\halfomega + a\right)} =
\overrightarrow{\CYB(\rho + b)}.$$ To finish the proof it is
enough to note that $[\pi, \pi] = p_*[\tilde\pi, \tilde\pi]$.
\hfill $\square$

\bigskip

\noindent {\it Proof of Theorem \ref{dyn}.} \ Let us construct the
inverse map. Suppose $(G/U, \pi)$ is a Poisson homogeneous $(G,
\pi_\rho)$-space. Set $b = \pi(\underline{e}) \in \wedge^2(\g /
\u) = \wedge^2\m$. The condition (c) implies that in fact $b \in
(\wedge^2(\g / \u))^\u = \mwmu$. Furthermore, (c) yields that
$\rho + b \in \halfomega + \mwmu$. By Proposition \ref{poiss}, we
have $\CYB(\rho + b) = 0$ in $\wedge^3(\g / \u)$, i.e., $\rho + b
\in \Momega$. Then, by Theorem \ref{ESmain}, there exists $r \in
\Dynr$ such that $r(0) = \rho + b$, and the image of $r$ in
$\Modr$ is uniquely determined. It is now easy to verify that $\pi
= \pi_r$. \hfill $\square$

\section{Twisting of Poisson homogeneous structures}\label{secttwist}

Assume again that $\g$ is an arbitrary finite-dimensional Lie
algebra over $\mathbb C$. Recall that a Lie bialgebra structure on
$\g$ is a 1-cocycle $\delta: \g \to \gwg$ which satisfies the
co-Jacobi identity. Denote by ${\mathcal D}(\g, \delta)$ the
classical double of $(\g, \delta)$.

We say that two Lie bialgebra structures $\delta_1$, $\delta_2$ on $\g$ are
in the same class if there exists a Lie algebra isomorphism $f: {\mathcal
D}(\g, \delta_1)\rightarrow {\mathcal D}(\g, \delta_2)$ which intertwines
the canonical forms $Q_i$ on ${\mathcal D}(\g, \delta_i)$, and such that the
following diagram is commutative:
%
%
%
$$ \xymatrix{ & \g \ar[dl]  \ar[dr] & \\
{\mathcal D}(\g,\delta_1) \ar[rr]^f & & {\mathcal D}(\g, \delta_2).}
$$
%
%
%
\begin{theorem} \label{twistth}
Two Lie bialgebra structures $\delta$, $\delta'$ on $\g$ are in
the same class if and only if $\delta' = \delta + ds$, where $s
\in \gwg$ and
\begin{equation} \label{twisteq}
\CYB(s) = \Alt(\delta \otimes \id)(s).
\end{equation}
\end{theorem}

\proof ($\Rightarrow$) Let us consider ${\mathcal D}(\g, \delta)$.
Then $\delta'$ is uniquely defined by a Lagrangian subalgebra $\l
\subset {\mathcal D}(\g, \delta)$ such that $\l \cap \g = 0$.
Clearly, $\l$ is the graph of a linear map $S: \g^\ast \to \g$.
Define an element $s = \sum_i s_i' \otimes s_i'' \in \g \otimes
\g$ via
\begin{equation}\label{S}
S(l)=\sum_i \langle l, s_i' \rangle s_i''
\end{equation}
for any $l \in \g^*$. Since $\l$ is Lagrangian, we see that $s$ is
skew-symmetric. Let us show that $\delta' = \delta + ds$.

Indeed, for any $a \in \g$, $l_1, l_2 \in \g^*$,
\[
\begin{split}
\langle \delta'(a), l_1 \otimes l_2 \rangle & = Q(\delta'(a),
(S(l_1) + l_1) \otimes (S(l_2) + l_2)) =\\ & = Q(a, [S(l_1) + l_1,
S(l_2) + l_2]) =\\ & = \langle a, [l_1, l_2] \rangle + Q(a,
[S(l_1), l_2]) + Q(a, [l_1, S(l_2)]),
\end{split}
\]
and
\[
\begin{split}
\langle a, [l_1, l_2] \rangle & = \langle \delta(a), l_1 \otimes
l_2 \rangle,\\
Q(a, [S(l_1), l_2]) & = \langle [a, S(l_1)], l_2 \rangle = \langle
[1 \otimes a, s], l_1 \otimes l_2 \rangle,\\
Q(a, [l_1, S(l_2)]) & = -\langle [a, S(l_2)], l_1 \rangle =
-\langle [1 \otimes a, s], l_2 \otimes l_1 \rangle =\\ & = \langle
[1 \otimes a, s^{21}], l_2 \otimes l_1 \rangle = \langle [a
\otimes 1, s], l_1 \otimes l_2 \rangle,
\end{split}
\]
where $\scpr$ is the canonical pairing between $\g$ and $\g^\ast$,
and $Q$ is the canonical bilinear form on ${\mathcal D}(\g,
\delta)$.

Now let $\{e_i\}$ be an arbitrary basis in $\g$ and $\{f^i\}$ be
its dual in $\g^\ast \subset {\mathcal D}(\g, \delta)$. Then the
canonical element $r_\delta=\sum_i e_i \otimes f^i \in {\mathcal
D}(\g, \delta)^{\otimes 2}$ satisfies the CYBE and $r_{\delta'} =
r_\delta + s$ satisfies the CYBE as well (since $r_\delta + s$ is
the canonical element for the double ${\mathcal D}(\g, \delta')$).
It is easy to show that these two facts imply (\ref{twisteq}).

($\Leftarrow$) $s \in \gwg$ defines $S: \g^\ast \rightarrow \g$
via (\ref{S}) and the graph of $S$ is $\l \subset {\mathcal D}(\g,
\delta)$, a Lagrangian subspace because $s$ is  skew-symmetric.
Let us prove that for any $l_1, l_2, l_3 \in \g^\ast$,
\[
\begin{split}
\langle l_1 \otimes l_2 \otimes l_3, \CYB(s) - \Alt(\delta \otimes
\id)(s) \rangle = Q([l_1+S(l_1), l_2+S(l_2)], l_3+S(l_3)).
\end{split}
\]
Let us verify that, for instance,
\[
\langle l_1\otimes l_2\otimes l_3, [s^{12}, s^{13}]\rangle =
Q([l_1, S(l_2)], S(l_3)).
\]
Indeed, if $s = \sum_i s_i' \otimes s_i''$, then we have
\[
[s^{12}, s^{13}] = \sum_{i, j} [s_i', s_j'] \otimes s_i'' \otimes
s_j''
\]
and
\[
\begin{split}
& \langle l_1\otimes l_2\otimes l_3, [s^{12}, s^{13}]\rangle =
\langle l_1,\sum_{i, j}[\langle s_i'', l_2\rangle s_i', \langle
s_j'', l_3\rangle s_j']\rangle = \\ & = \langle l_1,\sum_{i,
j}[\langle s_i', l_2\rangle s_i'', \langle s_j', l_3\rangle
s_j'']\rangle = \langle l_1, [S(l_2), S(l_3)]\rangle = \\ & =
Q(l_1, [S(l_2), S(l_3)]) = Q([l_1, S(l_2)], S(l_3)).
\end{split}
\]
Similarly,
\[
\begin{split}
- \langle l_1 \otimes l_2 \otimes l_3, (\delta \otimes
\id)(s)\rangle = - \langle [l_1,l_2]\otimes l_3, s\rangle =
\langle [l_1,l_2],S(l_3)\rangle = Q([l_1,l_2], S(l_3)),
\end{split}
\]
and so on. Since $Q([l_1,l_2],l_3)$ and $Q([S(l_1), S(l_2)],
S(l_3))$ vanish, the identity is proved.

Now it follows that $Q([l_1+S(l_1), l_2+S(l_2)], l_3+S(l_3))=0$
for any $l_1,l_2,l_3\in \g^{\ast}$. Since $\l$ is Lagrangian, we
conclude that $[l_1+S(l_1), l_2+S(l_2)]\in \l$ and hence $\l$ is a
subalgebra. Clearly, $\l$ defines $\delta' := \delta + ds$, and
this completes the proof of the theorem. \hfill $\square$

\begin{remark}
{\rm If we consider our Lie bialgebra $(\g, \delta)$ as a Lie
quasibialgebra, then $(\g, \delta + ds)$ is called ``twisting via
$s$''. The notions of Lie quasibialgebra and twisting via $s$ was
introduced by Drinfeld in \cite{Drinfeld QH}. The theorem above
can be also deduced from results of \cite{Drinfeld QH}.}
\end{remark}

Further, we are going to examine the effect of the twisting on Poisson
homogeneous spaces. First we recall some definitions and rather
well-known results.

Let $G$ be a connected complex Poisson-Lie group, $(\g,\delta)$
its Lie bialgebra, and ${\mathcal D}(\g) = {\mathcal
D}(\g,\delta)$ the corresponding classical double of $\g$ with the
canonical invariant form $Q$.

Recall that an action of $G$ on a Poisson manifold $M$ is called
Poisson if the defining map $G\times M\rightarrow M$ is a Poisson
map, where $G\times M$ is equipped with the product Poisson
structure. If the action is transitive, we say that $M$ is a {\it
Poisson homogeneous $G$-space}.

Let $M$ be a homogeneous $G$-space, and let $\pi$ be any bivector
field on $M$. For any $x \in M$ let us consider the map
$$\pi_x: T^\ast_xM\rightarrow T_xM,\ \pi_x(l)=(l \otimes
\id)(\pi(x)).$$
On the other hand, $M \cong G/H_x$ and $T_xM = \g/\h_x$, $T^\ast_xM =
(\g/\h_x)^\ast = \h^\perp_x \subset \g^\ast$, where $H_x$ is the stabilizer
of $x$, and $\h_x = \Lie H_x$. Therefore we can consider $\pi_x$ as a map
$\pi_x: \h^\perp_x \rightarrow \g/\h_x$.

Now let us consider the following set of subspaces in ${\mathcal
D}(\g) = \g \oplus \g^\ast$:
\begin{equation} \label{Lagrdef}
\l_x = \{ a + l \, | \, a \in \g, \, l \in \h^\perp_x,\, \pi_x(l)
= \overline{a} \},
\end{equation}
where $\overline{a}$ is the image of $a$ in $\g/\h_x$. Observe
that $\l_x$ are Lagrangian (i.e., maximal isotropic) subspaces,
and $\l_x \cap \g = \h_x$. The following result was obtained in
\cite{Drinfeld Lagr}.

\begin{theorem}[Drinfeld \cite{Drinfeld Lagr}] \label{Drinfeld}
$(M, \pi)$ is a Poisson homogeneous $G$-space if and only if for any $x \in
M$ the subspace $\l_x$ is a subalgebra of ${\mathcal D}(\g)$, and $\l_{gx} =
\Ad_g\l_x$ for all $g \in G$. \hfill $\square$
\end{theorem}


Now set $\delta' = \delta + ds$, where $s \in \gwg$ satisfies
(\ref{twisteq}). Then we have two Poisson-Lie groups, $(G,
\pi_\delta)$ and $(G, \pi_{\delta'})$, whose Lie bialgebras are
$(\g, \delta)$ and $(\g, \delta')$ respectively. Let $(M, \pi)$ be
a Poisson homogeneous $(G, \pi_\delta)$-space. Consider a bivector
field $\xi$ on $M$ defined by the formula $\xi(x) =$ the image of
$s$ in $\wedge^2 (\g / \h_x) = \wedge^2 T_x M$. Set $\pi' = \pi -
\xi$.

\begin{proposition} \label{homtwisted}
$(M, \pi')$ is a Poisson homogeneous $(G, \pi_{\delta'})$-space,
and thus one obtains a bijection between the sets of all Poisson
$(G, \pi_\delta)$- and $(G, \pi_{\delta'})$-homogeneous structures
on $M$.
\end{proposition}

\proof Theorem \ref{twistth} allows one to identify ${\mathcal
D}(\g, \delta)$ and ${\mathcal D}(\g, \delta')$. It is easy to
verify that under this identification the sets of Lagrangian
subspaces that correspond to $(M, \pi)$ and $(M, \pi')$ are the
same. This completes the proof, according to
Theorem~\ref{Drinfeld}. \hfill $\square$

\bigskip

Finally, we are going to generalize the main result of the
previous section to the twisted case. Assume that $(\g, \delta)$
is a quasitriangular Lie bialgebra, i.e., $\delta = d\rho$, where
$\rho \in \g \otimes \g$ and $\CYB(\rho) = 0$. It is easy to
verify that the condition (\ref{twisteq}) for an element $s \in
\gwg$ is equivalent to
\begin{equation} \label{twisteqtriang}
\CYB(s) + [\![\rho, s]\!] + [\![s, \rho]\!] = 0,
\end{equation}
where for $a, b \in \g^{\otimes 2}$ we set $[\![a, b]\!] = [a^{12}, b^{13}]
+ [a^{12}, b^{23}] + [a^{13}, b^{23}] \in \g^{\otimes 3}$ (i.e., $\CYB(a) =
[\![a, a]\!]$).

Fix $\Omega \in (S^2\g)^\g$ and assume that $\rho \in \halfomega +
\gwg$. As before, consider the Poisson-Lie group $(G,
\pi_\delta)$, where $\pi_\delta = \pi_\rho = \overrightarrow{\rho}
- \overleftarrow{\rho}$. Suppose $s \in \gwg$ satisfies
(\ref{twisteqtriang}). Set $\delta' = \delta + ds = d(\rho + s)$;
then $\pi_{\delta'} = \pi_{\rho, s} := \overrightarrow{\rho + s} -
\overleftarrow{\rho + s}$, and $(G, \pi_{\rho, s})$ is a
Poisson-Lie group.

Let $U$ be a connected Lie subgroup in $G$, and $\u = \Lie U$.
Consider $r \in \Dynr$. As usually, set $\tilde{\pi}_r =
\overrightarrow{r(0)} - \overleftarrow{\rho}$ and denote by
$\pi_r$ the natural projection of $\tilde{\pi}_r$ on $G/U$. By
Proposition \ref{Lu}, $(G/U, \pi_r)$ is a Poisson homogeneous $(G,
\pi_\rho)$-space. Set also $\tilde{\pi}_{r, s} = \tilde{\pi}_r -
\overleftarrow{s} = \overrightarrow{r(0)} - \overleftarrow{\rho +
s}$ and denote by $\pi_{r, s}$ its projection on $G/U$. According
to Proposition \ref{homtwisted}, $(G/U, \pi_{r, s})$ is a Poisson
homogeneous $(G, \pi_{\rho, s})$-space.

Moreover, if we combine Theorem \ref{dyn} and Proposition
\ref{homtwisted}, we get the following

\begin{theorem} \label{twistdyn}
Assume that $\u$, $\Omega$, and $\rho$ satisfy the conditions (a),
(b), and (c) from the previous section. Then the correspondence $r
\mapsto \pi_{r, s}$ is a bijection between $\Modr$ and the set of
all Poisson $(G, \pi_{\rho, s})$-homogeneous structures on $G/U$.
\hfill $\square$
\end{theorem}

Clearly, this can be reformulated as follows:

\begin{theorem} \label{twistdyn1}
Assume that $\u$ and $\Omega$ satisfy the conditions (a) and (b)
from the previous section. Suppose also that there exists $s \in
\gwg$ such that (\ref{twisteqtriang}) holds, and $\rho + s \in
\halfomega + \mwmu$. Then the correspondence $r \mapsto \pi_r$ is
a bijection between $\Modr$ and the set of all Poisson $(G,
\pi_\rho)$-homogeneous structures on $G/U$. \hfill $\square$
\end{theorem}

Let us apply our previous results to the triangular case.

\begin{corollary} \label{twistdyntriang}
Assume that $\u$ satisfies the condition (a) from the previous
section. Set $\Omega = 0$. Consider any $\rho \in \gwg$ that
satisfies the CYBE. Then the correspondence $r \mapsto \pi_r$ is a
bijection between $\Modr$ and the set of all Poisson $(G,
\pi_\rho)$-homogeneous structures on $G/U$.
\end{corollary}

\proof Set $s = -\rho$ and apply Theorem \ref{twistdyn1}. \hfill
$\square$

\section{Poisson homogeneous structures in the triangular case}\label{sectex}

Now assume that $\g$ is semisimple. Fix a Cartan subalgebra $\h \subset \g$
and denote by $\RRR$ the corresponding root system. In this section we apply
the results of the previous sections to the case $\u$ is a reductive Lie
subalgebra in $\g$ containing $\h$, $\Omega = 0$, and $\rho \in \gwg$ such
that $\CYB(\rho) = 0$.

To be more precise, consider $\UUU \subset \RRR$, and suppose $\u
= \h \oplus (\bigoplus_{\alpha \in \UUU} \g_\alpha)$ is a
reductive Lie subalgebra in $\g$. If this is the case, then we say
that a subset $\UUU \subset \RRR$ is {\it reductive} (i.e., $(\UUU
+ \UUU) \cap \RRR \subset \UUU$ and $-\UUU = \UUU$; see
\cite[Ch.~6, \S1.2]{Lie 3}). Condition (a) is satisfied since $\m
= \bigoplus_{\alpha \in \RRR \setminus \UUU} \g_\alpha$ is an
$\u$-invariant complement to $\u$ in $\g$.

Applying Corollary \ref{twistdyntriang} (and results of Etingof
and Schiffmann cited in Section \ref{sectdyn}), we get:

\begin{enumerate}
\item Any structure of a Poisson homogeneous $(G, \pi_\rho)$-space
on $G/U$ is of the form $p_*(\overrightarrow{x} -
\overleftarrow{\rho})$, where $x \in \Momega$.
\item If $x \in \Momega$, then there exists (a unique up to the
$\Mapszero$-action) $r\in \Dynr$ such that $r(0) = x$.
\end{enumerate}

Let us now describe $\Momega$ and thus get an explicit description
of all Poisson $(G, \pi_\rho)$-homogeneous structures on $G/U$.
Recall that in our case by definition $$\Momega =  \left\{ x \in
\mwmu \, \Bigm| \, \CYB(x) = 0 \ \mbox{in}\, \wedge^3(\g / \u)
\right\}.$$

We need to fix some notation. Fix a nondegenerate invariant
bilinear form ({\it invariant scalar product}) $\scpr$ on $\g$.
For any $\alpha \in \RRR$ choose $E_\alpha \in \g_\alpha$ such
that $\langle E_\alpha, E_{-\alpha} \rangle = 1$. Further, suppose
$\NNN$ is a reductive subset which contains $\UUU$. We say that $h
\in \h$ is {\it $(\NNN, \UUU)$-regular} if $\alpha (h) = 0$ for
all $\alpha \in \UUU$, and $\alpha (h) \neq 0$ for all $\alpha \in
\NNN \setminus \UUU$.

\begin{proposition} \label{Momega}
$x \in \Momega$ iff
\begin{equation} \label{xNh}
x = x_{\nh} = \sum_{\alpha \in \NNN \setminus \UUU}
\frac{1}{\alpha(h)}E_\alpha \otimes E_{-\alpha},
\end{equation}
where $\NNN$ is a reductive subset in $\RRR$ containing $\UUU$,
and $h \in \h$ is $(\NNN, \UUU)$-regular.
\end{proposition}

\proof First we calculate $\mwmu$. It is easy to see that $x \in
\mwm$ is $\h$-invariant iff it is of the form $$x = \sum_{\alpha
\in \RRR \setminus \UUU} x_\alpha \cdot E_\alpha \otimes
E_{-\alpha},$$ where $x_{-\alpha} = -x_\alpha$.

Define $c_{\alpha\beta}$ by the formula $[E_\alpha, E_\beta] =
c_{\alpha\beta}E_{\alpha + \beta}$ for $\alpha, \beta, \alpha +
\beta \in \RRR$.

Furthermore, suppose $\gamma \in \UUU$. One can easily verify
that the condition $\ad_{E_\gamma}(x) = 0$ is equivalent to the
following statement: for all $\alpha, \beta \in \RRR \setminus
\UUU$ such that $\alpha + \beta + \gamma = 0$ we have
$c_{\alpha\gamma}x_\alpha = c_{\beta\gamma}x_\beta$.

\begin{lemma} \label{c}
Suppose $\alpha, \beta, \gamma \in \RRR$, $\alpha + \beta + \gamma
= 0$. Then $c_{\alpha\gamma} + c_{\beta\gamma} = 0$.
\end{lemma}

\proof $c_{\alpha\gamma} = c_{\alpha\gamma}\langle E_{-\beta},
E_\beta \rangle = \langle [E_\alpha, E_\gamma], E_\beta \rangle =
\langle E_\alpha, [E_\gamma, E_\beta] \rangle =
-c_{\beta\gamma}\langle E_{\alpha}, E_{-\alpha} \rangle =
-c_{\beta\gamma}$. \hfill $\square$

\bigskip

Therefore we obtain

\begin{lemma} \label{l1}
$x \in \mwmu$ iff $$x = \sum_{\alpha \in \RRR \setminus \UUU}
x_\alpha \cdot E_\alpha \otimes E_{-\alpha},$$ where $x_{-\alpha}
= -x_\alpha$, and for all $\alpha, \beta \in \RRR \setminus \UUU$,
$\gamma \in \UUU$, $\alpha + \beta + \gamma = 0$, we have
$x_\alpha + x_\beta = 0$. \hfill $\square$
\end{lemma}

\begin{lemma} \label{l2}
Suppose $$x = \sum_{\alpha \in \RRR \setminus \UUU} x_\alpha \cdot
E_\alpha \otimes E_{-\alpha} \in \mwmu.$$ Then $x \in \Momega$ iff
the following condition holds: for all $\alpha, \beta, \gamma \in
\RRR \setminus \UUU$, $\alpha + \beta + \gamma = 0$, we have
$x_\alpha x_\beta + x_\beta x_\gamma + x_\gamma x_\alpha = 0$.
\end{lemma}

\proof One can check directly (using Lemma \ref{c}) that the
image of $\CYB(x)$ in $\wedge^3(\g/\u) = \wedge^3\m$ is equal to
$$\sum_{
\begin{array}{c}
\scriptstyle \alpha,\, \beta,\, \gamma \in \,\RRR \setminus \UUU,
\\ \scriptstyle \alpha + \beta + \gamma = 0
\end{array}}
c_{-\alpha, -\beta} \left(x_\alpha x_\beta + x_\beta x_\gamma +
x_\gamma x_\alpha\right)E_\alpha \otimes E_\beta \otimes
E_\gamma.$$ This immediately proves the lemma. \hfill $\square$

\bigskip

Now consider the following properties of the function $\RRR
\setminus \UUU \to \mathbb C$, $\alpha \mapsto x_\alpha$:

\begin{itemize}
\item[(d)] $x_{-\alpha} = -x_\alpha$ for all $\alpha \in \RRR
\setminus \UUU$.
\item[(e)] If $\alpha, \beta \in \RRR \setminus \UUU$,
$\gamma \in \UUU$, $\alpha + \beta + \gamma = 0$, then $x_\alpha +
x_\beta = 0$.
\item[(f)] If $\alpha, \beta, \gamma \in
\RRR \setminus \UUU$, $\alpha + \beta + \gamma = 0$, then
$x_\alpha x_\beta + x_\beta x_\gamma + x_\gamma x_\alpha = 0$.
\end{itemize}

\begin{lemma} \label{l3}
$x_\alpha$ satisfies (d)--(f) iff
\begin{equation}\label{x}
x_\alpha = \left\{
\begin{array}{cl}
1/\alpha(h), & \mbox{\rm if } \alpha \in \NNN \setminus \UUU
\\ 0, & \mbox{\rm if } \alpha \in \RRR \setminus \NNN ,
\end{array} \right.
\end{equation}
for a certain reductive subset $\NNN \subset \RRR$ such that $\NNN
\supset \UUU$, and $(\NNN, \UUU)$-regular element $h \in \h$.
\end{lemma}

\proof Suppose $x_\alpha$ satisfies (d)--(f). Set $\NNN = \UUU
\cup \{ \alpha \in \RRR \setminus \UUU \, | \, x_\alpha \neq 0
\}$. Let us prove that $\NNN$ is reductive. Using (d), we see that
$-\NNN = \NNN$. Further, suppose $\alpha, \beta \in \NNN$, $\gamma
\in \RRR$, $\alpha + \beta + \gamma = 0$. We have to verify that
$\gamma \in \NNN$. If $\alpha, \beta \in \UUU$, then also $\gamma
\in \UUU \subset \NNN$ (since $\UUU$ is reductive). If $\alpha \in
\UUU$, $\beta \in \NNN \setminus \UUU$, then $\gamma \in \RRR
\setminus \UUU$. Applying (e), we see that $x_\beta + x_\gamma =
0$. Since $x_\beta \neq 0$, we have $x_\gamma \neq 0$, i.e.,
$\gamma \in \NNN \setminus \UUU$. Finally, let $\alpha, \beta \in
\NNN \setminus \UUU$. Assume also that $\gamma \in \RRR \setminus
\UUU$ (we have nothing to prove in the case $\gamma \in \UUU$).
Using (f), we see that $x_\alpha \neq 0$, $x_\beta \neq 0$ imply
that $x_\gamma \neq 0$, i.e., $\gamma \in \NNN \setminus \UUU$.

Furthermore, set $y_\alpha = 1/x_\alpha$ for $\alpha \in \NNN
\setminus \UUU$. Suppose $\alpha, \beta, \gamma \in \NNN \setminus
\UUU$, $\alpha + \beta + \gamma = 0$. Then $y_\alpha + y_\beta +
y_\gamma = 0$ according to (f). This means that $y_\alpha = \alpha
(h)$ for some $h \in \h$.

Finally, we prove that $h$ is $(\NNN, \UUU)$-regular. By
construction, $\alpha(h)\neq 0$ for all $\alpha \in \NNN \setminus
\UUU$. Now assume that $\gamma \in \UUU$. Take any $\alpha \in
\NNN \setminus \UUU$ (note that if $\NNN = \UUU$, then we have
nothing to prove here), and set $\beta = -(\alpha + \gamma)$.
Obviously, $\beta \in \NNN \setminus \UUU$. By (e), we have $0 =
x_\alpha + x_\beta = 1/\alpha(h) + 1/\beta(h)$, i.e., $\gamma(h) =
0$.

Conversely, if $x_\alpha$ is of the form (\ref{x}), then the
conditions (e)--(f) can be verified without difficulties. \hfill
$\square$

\bigskip

The last lemma proves the proposition. \hfill $\square$

\begin{remark}
{\rm We note that Lemmas \ref{l1}, \ref{l2}, and \ref{l3} are
essentially contained in \cite{DGS}.

In \cite{DGS}, among other results, the {\it symplectic}
$G$-invariant structures on $G/U$ are classified if $U$ is a Levi
subgroup of $G$. Actually, in this case there exists a
$G$-equivariant symplectomorphism from $G/U$ to a semisimple
coadjoint $G$-orbit equipped with the Kirillov-Kostant-Souriau
bracket.

Moreover, it is easy to show that if $G/U$ has a $G$-invariant
symplectic structure, then $U$ is a Levi subgroup. Indeed, let
$p_*\overrightarrow{x_{\nh}}$ (where $x_{\nh}$ is defined by
(\ref{xNh})) be a $G$-invariant Poisson structure on $G/U$.
Obviously, it is symplectic iff $\NNN = \RRR$. Since $h$ is
$(\RRR, \UUU)$-regular, i.e., $\alpha(h) = 0$ for all $\alpha \in
\Span \UUU$ and $\alpha(h) \neq 0$ for all $\alpha \in \RRR
\setminus \UUU$, we see that $(\Span \UUU) \cap \RRR = \UUU$. It
is well known that the latter condition is equivalent to the fact
that $U$ is a Levi subgroup.

Let us also remark that the existence of reductive non-Levi
subgroups is the main difference between the triangular and the
strictly quasitriangular cases. Indeed, suppose $U$ is a Cartan
subgroup. Then in the triangular case the Poisson homogeneous
structures on $G/U$ relate to {\it all}\/ reductive subgroups of
$G$, while in the strictly quasitriangular case they relate to the
Levi subgroups only (see \cite{Kar1, Lu CDYBE}).}
\end{remark}

Now we are going to describe the Lagrangian subalgebras
corresponding to the Poisson $(G, \pi_\rho)$-homogeneous
structures on $G/U$. Since the Lie bialgebras
corresponding to $(G, \pi_\rho)$ are all in the same class, we may
assume without loss of generality that $\rho = 0$. It is clear
that the double of our Lie bialgebras is $\geps = \g \oplus
\mathfrak g\varepsilon$, where $\varepsilon^2 = 0$ (see Appendix A
for details).

Suppose $\rho = 0$. Assume that $\NNN$ and $h$ are as in
Proposition \ref{Momega}. Set $\pi_{\nh} =
p_*\overrightarrow{x_{\nh}}$, where $x_{\nh}$ is defined by
(\ref{xNh}). By $\l_{\nh}$ denote the Lagrangian subalgebra
corresponding to $(G/U, \pi_{\nh})$ at the base point
$\underline{e}$.

\begin{proposition}
$\l_{\nh} = \u \oplus \left( \bigoplus_{\alpha \in \RRR \setminus
\NNN} \varepsilon\g_\alpha \right) \oplus \left( \bigoplus_{\alpha
\in \NNN \setminus \UUU} (1 - \alpha (h) \varepsilon)\mathfrak
g_\alpha \right)$ (cf.~Proposition \ref{lnb} below).
\end{proposition}

\proof By definition (see (\ref{Lagrdef})), $$\l_{\nh} = \{ a +
b\varepsilon \, | \, a \in \g, b \in \u^\perp = \m, (b \otimes
1)(x_{\nh}) = \overline{a} \}, $$ where $\overline{a}$ is the
image of $a$ in $\g/\u = \m$. Suppose $b = E_\alpha$, where
$\alpha \in \RRR \setminus \UUU$. Then $$(b \otimes 1)(x_{\nh}) =
\left\{
\begin{array}{cl}
-\frac{1}{\alpha(h)}E_\alpha, & \mbox{\rm if } \alpha \in \NNN
\setminus \UUU \\ 0, & \mbox{\rm if } \alpha \in \RRR \setminus
\NNN.
\end{array} \right. $$
This completes the proof. \hfill $\square$

\appendix

\section{Lagrangian subalgebras in $\geps$}

Let $\g$ be a semisimple complex Lie algebra, $G$ a connected Lie
group such that $\Lie G = \g$. Fix an invariant scalar product
$\scpr$ on $\g$. Consider the complex Lie algebra $\geps = \g
\otimes_{\mathbb C} \mathbb C[\varepsilon] = \g \oplus \mathfrak
g\varepsilon$, where $\mathbb C[\varepsilon] = \mathbb C[x]/(x^2)$
is the algebra of dual numbers. We identify $\g$ with $\g \otimes
1 \subset \geps$. Equip $\geps$ with the invariant scalar product
defined by $$\langle a + b\varepsilon, c + d\varepsilon \rangle =
\langle a, d \rangle + \langle b, c \rangle.$$ Then the pair
$(\geps, \g)$ is a Manin pair.

Recall that a Lie subalgebra $\l \subset \geps$ is called {\it
Lagrangian} if it is a maximal isotropic subspace in $\geps$.

Let $\n \subset \g$ be a Lie subalgebra, $B$ be a $\mathbb
C$-valued 2-cocycle on $\n$. By $\Pairs$ denote the set of all
such pairs $(\n, B)$. Define $f: \n \to \n^*$ by
\begin{equation} \label {Bf}
\langle f(x), y \rangle = B(x,y)
\end{equation}
(here we identify $\n^*$ with $\g/\n^\perp$ via $\scpr$). One can
easily see that $f$ is a skew-symmetric 1-cocycle (with respect to
the coadjoint action of $\n$ on $\n^*$). Set $$\l (\n, B) = \{ a +
b\varepsilon \, | \, a \in \n, b \in \g, f(a) = \bar{b} \} \subset
\geps, $$ where $\bar{b}$ is the image of $b$ in $\n^* =
\g/\n^\perp$. Obviously, $\u := \Ker B = \Ker f$ is a Lie
subalgebra in $\n$.

We denote by $\Lagr$ the set of all Lagrangian subalgebras in
$\g$. Note that $G$ acts naturally on $\Pairs$ and on $\Lagr$.

\begin{theorem} \label{Pairs}

1. $\l (\n, B)$ is a Lagrangian subalgebra in $\geps$ and $\l
(\n, B) \cap \g = \u$.

2. The mapping $(\n, B) \mapsto \l (\n, B)$ is
a $G$-equivariant bijection between $\Pairs$ and $\Lagr$.
\end{theorem}

\proof Suppose $\l \in \Lagr$. Denote by $\n$ the projection of
$\l$ onto $\g$ along $\g\varepsilon$. Then $\l \subset \n \oplus
\g\varepsilon$ and $\l \cap (\g\varepsilon) = \n^\perp\varepsilon$.
Consider
$$\bar{\l} := \l / (\n^\perp\varepsilon) \subset (\n \oplus
\g\varepsilon) / (\n^\perp\varepsilon) = \n \oplus
\n^*\varepsilon.$$
Since $\dim \bar{\l} = \dim \l - \dim \n^\perp
= \dim \g - (\dim \g - \dim \n) = \dim \n$, we see that $\bar{\l}$
is the graph of a linear map $f: \n \to \n^*$, i.e., $$\bar{\l} =
\{ a + f(a)\varepsilon \, | \, a \in \n \}.$$ This yields that
$$\l = \{ a + b\varepsilon \, | \, a \in \n, b \in \g, f(a) =
\bar{b} \} \subset \geps, $$ where $\bar{b}$ is the image of $b$
in $\n^* = \g/\n^\perp$.

Now let $a + b\varepsilon, c + d\varepsilon \in \l$ (i.e., $a, c
\in \n$, $f(a) = \bar{b}$, $f(c) = \bar{d}$). Since $\l$ is a Lie
subalgebra, we have
$$\l \ni [a + b\varepsilon, c + d\varepsilon] = [a, c] + ([a, d] +
[b, c])\varepsilon.$$
Therefore
$$f([a, c]) = \overline{[a, d] + [b, c]} = [a, f(c)] + [f(a),
c],$$
i.e., $f$ is a 1-cocycle. Since $\l$ is isotropic, we have
$$0 = \langle a + b\varepsilon, c + d\varepsilon \rangle = \langle
a, d \rangle + \langle b, c \rangle = \langle a, f(c) \rangle +
\langle f(a), c \rangle,$$
i.e., $f$ is skew-symmetric. Finally,
define $B$ by (\ref{Bf}). It is easy to check that $B$ is a
2-cocycle.

Conversely, $\l (\n, B)$ is a Lie subalgebra since $\n$ is a Lie
subalgebra and $f$ is a 1-cocycle (recall that $f$ and $B$ are
connected via (\ref{Bf})); $\l (\n, B)$ is isotropic since $f$ is
skew-symmetric; finally, $\l (\n, B)$ is Lagrangian since $\dim \l
(\n, B) = \dim \n + \dim \n^* = \dim \g$.

The fact that $\l (\n, B) \cap \g = \u$ is obvious. \hfill
$\square$

\bigskip

Now fix a Lie subalgebra $\u \subset \g$. Set $$\Pairsu = \{ (\n,
B) \in \Pairs \, | \, \n \supset \u, \Ker B = \u \}, $$ $$\Lagru =
\{ \l \in \Lagr \, | \, \l \cap \g = \u \}.$$ Denote by $N(\u)$
the normalizer of $\u$ in $G$. Clearly, $N(\u)$ acts on $\Pairsu$
and $\Lagru$.

\begin{corollary} \label{Pairsu}
The mapping $(\n, B) \mapsto \l (\n, B)$ is a $N(\u)$-equivariant
bijection between $\Pairsu$ and $\Lagru$. \hfill $\square$
\end{corollary}



As before, fix a Cartan subalgebra $\h \subset \g$, and denote by
$\RRR$ the corresponding root system. Consider a reductive subset
$\UUU \subset \RRR$ and set $\u = \h \oplus (\bigoplus_{\alpha \in
\UUU} \g_\alpha)$. We would like to describe more explicitly the
set $\Lagru$. By Corollary \ref{Pairsu}, it
is sufficient to describe the set $\Pairsu$.


\begin{theorem} \label{Descru}
Suppose $(\n, B) \in \Pairs$. Then $(\n, B) \in \Pairsu$ if and
only if $\n = \h \oplus (\bigoplus_{\alpha \in \NNN} \g_\alpha)$
is a reductive Lie subalgebra in $\g$ that contains $\u$, and
$B(x, y) = \langle h, [x, y] \rangle$ (i.e., $B$ is a
2-coboundary), where $h \in \h$ is $(\NNN, \UUU)$-regular (for
the definition see the previous section).
\end{theorem}

\begin{remark}
{\rm Suppose $h \in \h$. It is clear that $B(x, y) = \langle h,
[x, y] \rangle$ depends only on the image of $h$ in $\h / \z
(\n)$, where $\z (\mathfrak n)$ is the center of $\n$. Note also
that $\z (\mathfrak n) = \{ h \in \h \, | \, \alpha (h) = 0 \ {\rm
for \ all}\ \alpha \in \NNN \}$.}
\end{remark}

\noindent {\it Proof of Theorem \ref{Descru}.} \ Suppose $(\n, B)
\in \Pairsu$, i.e., $\n \supset \u$, $\Ker B = \u$. Since $\n
\supset \h$, we see that $\n = \h \oplus (\bigoplus_{\alpha \in
\NNN} \g_\alpha)$ for some $\NNN \subset \RRR$. Clearly, $\UUU
\subset \NNN$.

\begin{lemma} \label{B}
If $\alpha, \beta \in \NNN$, $\alpha + \beta \neq 0$, then $B(x,
y) = 0$ for all $x \in \g_\alpha$, $y \in \g_\beta$.
\end{lemma}

\proof If $x \in \g_\alpha$, $y \in \g_\beta$, $h \in \mathfrak
h$, then, since $B$ is a 2-cocycle, we have $B([x, y], h) + B([y,
h], x) + B([h, x], y) = 0$, i.e., $B(h, [x, y]) = (\alpha +
\beta)(h) \cdot B(x, y)$. Since $\Ker B \supset \h$ and $\alpha +
\beta \neq 0$, we see that $B(x, y) = 0$. \hfill $\square$

\bigskip

Now we continue the proof of the theorem. If $\alpha \in \NNN$, but
$-\alpha \notin \NNN$, then, by Lemma \ref{B}, we see that
$\g_\alpha \subset \Ker B=\u$. Then $\pm\alpha \in \UUU\subset\NNN$
because $\UUU$ is reductive. This contradiction
proves that $-\NNN = \NNN$, i.e., $\n$ is
reductive.

Let us prove that $B$ is a 2-coboundary. Recall that $H^2(\n,
\mathbb C) = \wedge^2 \z (\n)$ (see \cite{Lie 2}). To be more
precise, any 2-cocycle $B$ can be presented uniquely in the form
$B' + B''$, where $B'$ is a 2-coboundary, and $B''(x, y) = \langle
u, x \otimes y \rangle$ for $u \in \wedge^2 \z (\n)$. Assume that
$B'' \neq 0$. Then there exists $a \in \z (\n) \subset \h$ such
that $a \notin \Ker B''$. Since $B'$ is a 2-coboundary, we see
that $a \in \Ker B'$. Therefore $a \notin \Ker B$, and we get a
contradiction. This means that $B (x, y) = \langle h, [x, y]
\rangle$, where $h \in \n$.

It remains to prove that $h$ is a $(\NNN, \UUU)$-regular element
of $\h$. Suppose $\alpha \in \NNN$, $x \in \g_\alpha$, $h' \in
\h$. Since $\Ker B \supset \h$, we have $0 = B (h', x) = \langle
h, [h', x] \rangle = \alpha (h') \cdot \langle h, x \rangle$.
Therefore $h$ is orthogonal to $\g_\alpha$ for all $\alpha \in
\NNN$. This implies that $h \in \h$.

If $\alpha \in \NNN$, $x \in \g_\alpha$, $y \in \n$, then $B(x, y)
= \langle [h, x], y \rangle = \alpha (h) \cdot \langle x, y
\rangle$. This shows that $\g_\alpha \subset \Ker B$ iff $\alpha
(h) = 0$. Therefore $\Ker B = \u$ iff $h$ is $(\NNN,
\UUU)$-regular.

The converse statement of the theorem can be verified
directly. \hfill $\square$

\bigskip

Suppose $\NNN$ is a reductive subset in $\RRR$ containing $\UUU$.
By $\n$ denote the reductive Lie subalgebra in $\g$ that
corresponds to $\NNN$. Consider a $(\NNN, \UUU)$-regular element
$h \in \h$. Denote by $B$ the 2-coboundary which corresponds to
$h$ (see Theorem \ref{Descru}).

\begin{proposition} \label{lnb}
$\l (\n, B) = \h \oplus \left( \bigoplus_{\alpha \in \RRR
\setminus \NNN} \varepsilon\mathfrak g_\alpha \right) \oplus
\Bigl( \bigoplus_{\alpha \in \NNN} (1 + \alpha (h)
\varepsilon)\g_\alpha \Bigr) = \u \oplus \left( \bigoplus_{\alpha
\in \RRR \setminus \NNN} \varepsilon\g_\alpha \right) \oplus
\left( \bigoplus_{\alpha \in \NNN \setminus \UUU} (1 + \alpha (h)
\varepsilon)\mathfrak g_\alpha \right)$.
\end{proposition}

\proof Direct calculations. \hfill $\square$


\section{From dynamical to constant r-matrices}

In this appendix we describe a procedure that leads from dynamical
to constant r-matrices.

Let $\g$ be any complex finite-dimensional Lie algebra, $\u \subset \g$
a Lie subalgebra. Assume that there exists a Lie subalgebra
$\v \subset \g$ such that $\g = \u \oplus \v$ as a vector space.
Let $r$ be a classical dynamical r-matrix for the pair $(\g, \u)$.
Set $\Omega := r + r^{21}$ and assume that $\Omega \in (\u \otimes \u)
\oplus (\v \otimes \v)$. By $v$ denote the image of $r(0)$ under
the projection onto $\v \otimes \v$ along
$\g \otimes \u + \u \otimes \g$.

\begin{proposition}\label{d-c}
$\CYB(v) = 0$.
\end{proposition}

\proof Let $r(0) = \sum_i r_i' \otimes r_i'' \in \g \otimes
\g$. Decompose $r_i' = a_i' + b_i'$, $r_i'' = a_i'' + b_i''$, where
$a_i', a_i'' \in \u$, $b_i', b_i'' \in \v$. Then $v = \sum_i b_i'
\otimes b_i'' \in \v \otimes \v$.

Since $r$ satisfies the CDYBE, we see that $\CYB(r(0)) \in
\g \otimes \g \otimes \u + \g \otimes \u \otimes \g +
\u \otimes \g \otimes \g$, i.e., the $\v \otimes \v \otimes \v$-component
of $\CYB(r(0))$ is zero. On the other hand, one can calculate directly that
the $\v \otimes \v \otimes \v$-component of $\CYB(r(0))$ equals
$\CYB(v) + \sum_{k = 1}^6 A_k$, where
$$A_1 = \sum_{i, j} [a_i', b_j']_\v \otimes b_i'' \otimes b_j'',$$
$$A_2 = \sum_{i, j} b_i' \otimes [a_i'', b_j']_\v \otimes b_j'',$$
$$A_3 = \sum_{i, j} b_i' \otimes b_j' \otimes [a_i'', b_j'']_\v,$$
$$A_4 = \sum_{i, j} [b_i', a_j']_\v \otimes b_i'' \otimes b_j'',$$
$$A_5 = \sum_{i, j} b_i' \otimes [b_i'', a_j']_\v \otimes b_j'',$$
$$A_6 = \sum_{i, j} b_i' \otimes b_j' \otimes [b_i'', a_j'']_\v;$$
here we denote by $x_\v$ the $\v$-component of $x \in \g$.
We have to show that $\sum_{k = 1}^6 A_k = 0$.

Since $r$ is $\u$-equivariant, we see that $r(0)$ is $\u$-invariant.
Now consider any $a \in \u$. Taking the $\v \otimes \v$-component of
$\ad_a r(0)$, we have
\begin{equation} \label{uinv}
\sum_{i} [a, b_i']_\v \otimes b_i'' + \sum_{i} b_i' \otimes [a, b_i'']_\v = 0.
\end{equation}
This implies immediately that $A_2 + A_3 = 0$, $A_4 + A_5 = 0$.

Now notice that the condition $\Omega \in (\u \otimes \u)
\oplus (\v \otimes \v)$ means that $\sum_i a_i' \otimes b_i'' =
-\sum_i a_i'' \otimes b_i'$. This yields that
\[
\begin{split}
A_1 & = -\sum_{i, j} [a_i'', b_j']_\v \otimes b_i' \otimes b_j'' =\\
& = \sum_{i, j} [b_i', a_j'']_\v \otimes b_j' \otimes b_i''.
\end{split}
\]
Combining this with (\ref{uinv}), we see that $A_1 + A_6 = 0$.
\hfill $\square$

\begin{remark}
{\rm Suppose $g \in \Mapszero$, and let $r$ be a classical dynamical
r-matrix for the pair $(\g, \u)$. Since $r^g(0) - r(0) \in \u \otimes
\g + \g \otimes \u$, we see that the $\v \otimes \v$-components
of $r(0)$ and $r^g(0)$ are the same. Therefore under
conditions of Proposition \ref{d-c} we get a map from $\Modr$
to the set of all classical r-matrices for $\v$.}
\end{remark}

\begin{example}
{\rm Let $\g = \mathfrak{sl}(n)$.
Suppose $\h$ is the standard Cartan subalgebra
in $\g$, i.e., the subalgebra of diagonal matrices. Consider a parabolic
subalgebra $\p \subset \g$ of the form
$$
\p=\left\{\left(\begin{array}{cccc}
*&\ldots&*&0\\
\vdots&\ddots&\vdots&\vdots\\
*&\ldots&*&0\\
*&\ldots&*&*
\end{array}\right)\right\}.
$$
Let $g = \id - E_{1n} - E_{2n} - \ldots - E_{n-1, n} \in \mbox{GL}(n)$, where
$E_{ij}$ is the matrix that has $1$ at $(i, j)$-entry and $0$ elsewhere.
It is easy to show that $\g = \h \oplus g\p g^{-1}$ as a vector space
(see \cite{GS}). Let us apply Proposition \ref{d-c} to this situation.

Fix $h_1, \ldots, h_n \in \mathbb C$ such that $h_i \neq h_j$ for all
$i \neq j$ and $h_1 + \ldots + h_n = 0$. It follows from the results of
Etingof and Varchenko \cite{E-V} (see also Section 4) that there exists
a classical dynamical r-matrix $r$ for the pair
$(\g, \h)$ such that
$$r(0) = \sum_{1 \leq i < j \leq n}\frac{1}{h_i - h_j}E_{ij}\wedge E_{ji}.$$
By $p$ denote the projection of $\g \otimes \g$ onto
$g(\p \otimes \p)g^{-1}$ along $\g \otimes \h + \h \otimes \g$.
Since $r$ is skew-symmetric, the conditions of Proposition \ref{d-c} are
obviously satisfied. Therefore $v := g^{-1} p(r(0)) g$ is a triangular
classical r-matrix for $\p$. It is not hard to calculate that
\[
\begin{split}
v = & \sum_{1 \leq i < j \leq n-1}\frac{1}{h_i - h_j}
(E_{ij} - D_i)\wedge (E_{ji} - D_j) +\\
+ & \sum_{i = 1}^{n-1}\frac{1}{h_i - h_n}D_i\wedge
\left(\sum_{1 \leq k \leq n,\ k \neq i}E_{ki}\right);
\end{split}
\]
here $D_i = \mbox{diag}\left(\frac{1}{n}, \ldots, \frac{1}{n}, -\frac{n-1}{n},
\frac{1}{n}, \ldots, \frac{1}{n}\right) \in \h$
(with $-\frac{n-1}{n}$ at $i$th place).}
\end{example}

\end{sloppy}

\end{document}